\documentclass [12pt]{article}

\usepackage{latexsym,amsthm,amsfonts,amsmath}

\newtheorem{theo}{Theorem}[section]

\newtheorem{prop}[theo]{Proposition}

\newtheorem{coro}[theo]{Corollary}

\newcommand\N{\ensuremath{\mathbb{N}}}

\newcommand\Z{\ensuremath{\mathbb{Z}}}

\renewcommand\P{\ensuremath{\mathbb{P}}}

\title{Oriented Random Walk on the Heisenberg Group and Percolation}
\author{Itai Benjamini and Oded Schramm}

\begin{document}
\maketitle

\begin{abstract}
It is shown that oriented random walk on the Heisenberg group
admits exponential intersection tail. As a corollary we get that on any
transitive graph of polynomial volume growth, which is
not a finite extension of $\Z^1, \Z^2$, the infinite
cluster of percolation with retention parameter $p$, close enough to
$1$, is transient.
\end{abstract}

\section{Introduction}
The study of random walks on discrete groups is rather advanced
(see for instance Hebisch and Saloff-Coste (1993)).  Recently a
study of percolation and other related models on groups was
initiated (see Benjamini and Schramm (1996) for an introduction).
For the theory of percolation on groups to mature it is useful
to have some interesting specific groups in which percolation is
analyzed. Since the discrete Heisenberg group is the smallest (in
the sense of volume growth) infinite non commutative discrete group, which
is not a finite extension of $\Z^d$ $(d \leq 4)$, it is a
natural example to consider.  We will see below that the analysis of
oriented random walks on the Heisenberg group is useful for the study
of percolation on transitive graphs of polynomial volume growth.

\noindent
{\bf Description of the Heisenberg group and its Cayley graph.}
The Heisenberg group $H$ has the presentation
$H=\langle a,b,c : c=[a,b],\,[a,c]=[b,c]=e\rangle$,
where $e$ denotes the identity element.
We now describe the Cayley graph $G_H$ of $H$ with respect to
the generators $a$ and $b$.
It is easy to verify that each element
$h\in H$ has a unique representation of the form
$h=a^nb^m c^k$, where $n,m,k\in\Z$.  Thus $h$
can be identified with the point $(n,m,k)\in\Z^3$,
and the vertices of $G_H$ may be taken to be $\Z^3$.
The edges then are
$$
E =
\bigl\{ [(x,y,z), (x ,y+1,z)] : x,y,z\in\Z \bigr\} \cup
\bigl\{ [(x,y,z), (x+1 ,y,z- y)] : x,y,z\in\Z \bigr\}
\,.
$$
With this set of directed edges $G_H$ may be viewed as a {\bf directed}
graph.
A {\bf directed} path on $G_H$ is a path that respects the orientation of
its edges.

Since $G_H$ contains a copy of $\Z^2$, $p_c(G_H) <1$.
Several facts regarding percolation
on $\Z^d$, such as uniqueness of the infinite cluster or exponential decay
for the connectivity function below the critical probability, holds for $G_H$
with similar proofs. It is of interest to prove an analogue of the
Grimmett-Marstrand theorem regarding percolation in a slab for $G_H$.
As this is the missing part in applying renormalization arguments for
percolation on $G_H$.
(for background on percolation see
Grimmett (1989)). Grimmett Kesten and Zhang (1993) proved transience of
the super critical infinite cluster in $\Z^d$, $d \geq 3$.
Benjamini, Pemantle and Peres (1998) gave a very different proof
using a measure on paths such that the number of intersections of two
paths chosen independently according to the measure, has an exponential tail.

In this note it is shown that the
probability that two independent oriented random
walks on $G_H$  will intersect $n$ times decays exponentially in $n$.
As a corollary we get that $G_H$ is transient and even the super critical
infinite percolation cluster, for $p$ close enough to $1$, is  a.s. transient.

We start with some definitions and a result from  Benjamini, Pemantle and
Peres (1998) and then formulate the theorem.
\medskip

\noindent{\bf Definitions.}
\vspace{-.1in}
\begin{enumerate}
\item
Let $G=(V_G,E_G)$ be an infinite directed graph with all vertices of
finite degree and let $v_0 \in V_G$. Denote by
$\Upsilon= \Upsilon(G,v_0)$ the collection of infinite directed paths in $G$ which
emanate from $v_0$ and tend to infinity (i.e., the paths in $\Upsilon$ visit any
vertex at most finitely many times).
 The set $\Upsilon(G,v_0)$, viewed as a subset of
$E_G^\Z$, is a Borel set in the product topology.
\item
Let $0 < \theta <1$. A Borel probability measure $\mu$ on $\Upsilon(G,v_0)$
has   {\bf Exponential intersection tails} with parameter $\theta$
(in short, EIT($\theta$)) if there exists $C$ such that

$$ 
\mu \times \mu \Big\{(\varphi,\psi): | \varphi \cap \psi | \geq n \Big\}
\leq C \theta^n
$$
for all $n$, where  $| \varphi \cap \psi |$ is the number of edges in the
intersection of $\varphi$ and $\psi$.
\item
If such a measure $\mu$ exists for some  basepoint $v_0$ and some $\theta<1$, then
we say that $G$ {\em admits random paths with\/} EIT($\theta$).
Analogous definitions apply to undirected graphs.
\item {\bf Oriented percolation}  with parameter $p \in (0,1)$ on the directed
graph $G$ is the process where each edge of $G$  is independently declared
 {\em open} with probability $p$ and {\em closed} with probability $1-p$.
The union of all directed open paths emanating from $v$ will be called the
{\bf oriented open cluster} of $v$ and denoted $C(v)$.
\item
A subgraph $\Lambda$ of $G$ is called {\bf transient}
if when the orientations on the edges are ignored,
$\Lambda$ is connected and simple random walk on it is a transient Markov
chain. As explained in Doyle and Snell (1984), the latter property
is equivalent to
finiteness of the effective resistance from a vertex of $\Lambda$
to infinity, when each edge of $\Lambda$ is endowed with a unit resistor.
\end{enumerate}

The following proposition is from Benjamini, Pemantle and Peres (1998).

\begin{prop} \label{prop:trans}
Suppose a directed graph $G$ admits random paths with EIT($\theta$).
Consider oriented percolation on $G$ with parameter $p$.
If $p>\theta$ then with probability 1 there is a vertex $v$ in $G$
such that the directed open cluster $C(v)$ is transient.
\end{prop}

Recall that a path $\{\Gamma_n\}$ in $\Z^d$ is called {\em oriented} if each
 increment $\Gamma_{n+1}-\Gamma_n$ is one of the $d$ standard basis vectors.
 The difference of two independent, uniformly chosen, oriented paths
in $\Z^d$ is a random walk with increments generating the $d-1$
dimensional hyperplane
$\{\sum_{i=1}^d x_i =0 \}$. For $d \geq 4$, this random walk is transient;
let $\theta_d<1$ denote  its return probability to the origin. As noted by
Cox and Durrett (1983), it follows that
 the uniform measure on oriented paths in $\Z^d$ has EIT($\theta_d$).
(They attribute the idea of applying this in percolation to  H. Kesten.)

\begin{theo}
\label{main}
The uniform measure on oriented paths in $G_H$ has
exponential intersection tail.
\end{theo}
\medskip

Where the uniform measure on oriented paths is the infinite symmetric product
measure on $\{a,b\}^{\N}$, were each word corresponds to an oriented path
naturally.

It is well known (to group theorists) that any
finitely generated infinite
nilpotent group, which is not a finite extension of $\Z$ or $\Z^2$,
contains either the discrete Heisenberg group or $\Z^3$ as subgroups.
(For finitely generated nilpotent groups  the torsion group is finite,
assuming this, it is an exercise to show the above, see Kurosh (1956)).
By a Trofimov's generalization of Gromov's celebrated theorem
(Trofimov (1986)), any
transitive graph of polynomial volume growth  is a nilpotent cover of
a finite graph. Thus in particular contains a graph which is rough-
isometric to either $G_H$ or $\Z^3$. Thus we get

\begin{coro}
\label{perc}
The infinite cluster for percolation with retention parameter $p$
on a transient transitive graph of polynomial volume growth is transient for
$p$ close enough to $1$.
\end{coro}
\medskip

Benjamini, Lyons and Schramm (1998) conjectured that a super
critical infinite percolation cluster on any transitive transient graph
is transient. The conjecture was verified there for nonamenable Cayley
graphs.

\section{Proof and further remarks}

Consider two oriented random walks
$\alpha=(\alpha(0),\alpha(1)\dots)$ and $\beta=(\beta(0),\beta(1)\dots)$.
Let $\alpha_j=0$ if the $j$'th edge of $\alpha$ is of the
form $[p,pa]$, and $\alpha_j=1$ if it is of the form
$[p,pb]$.  Similarly define $\beta_j$.
If $v=(v_1,v_2,v_3)\in\Z^3$ is the location of $\alpha$
after $k$ steps, then $v_1+v_2=k$.
Consequently, intersections of $\alpha$ and $\beta$
correspond to $k$'s such that $\alpha(k)=\beta(k)$.

For $\alpha(k)=\beta(k)$,
we must have $\sum_{j<k}\alpha_j= \sum_{j<k}\beta_j$
and also $\sum_{j<k}\sum_{i<j} (1-\alpha_i)\alpha_j$
must be equal to the corresponding expression with $\beta$.
$$
\sum_{j<k}\sum_{i<j} (1-\alpha_i)\alpha_j= \sum_{j<k} j \alpha_j -\sum_j\sum_{i<j}\alpha_i \alpha_j .
$$ 

Since $\sum_j\sum_{i<j}\alpha_i \alpha_j$ depends
only on the number of $i$'s for which $\alpha_i =1$,
$\alpha(k)=\beta(k)$ convertible to
$\sum_{j<k} j\alpha_j= \sum_{j<k} j\beta_j$ together with
$\sum_{j<k}\alpha_j= \sum_{j<k}\beta_j$.
The probability for the latter is of order $k^{-1/2}$.
Given the latter, the probability of the former is of
order $k^{-3/2}$.  However, the following easier estimate
suffices:
$$
\P\left[\sum_{j<k} j\alpha_j= \sum_{j<k} j\beta_j\right] \leq {1\over k}
\,.
$$
Note that
$$
\P\left[\sum_{j<k} j\alpha_j= \sum_{j<k} j\beta_j\right]
\le \max_n
\P\left[\sum_{j<k} j\alpha_j=n\right]
$$
Hence the above estimate can be obtained by observing that
$\sum_{j<\log_2 k} 2^j \alpha_{2^j}$ is uniformly distributed
in an interval (of $\Z$) of size at least $(k/2)-1$.
Together with the other condition, this gives
a probability of $O(k^{-3/2})$ for an intersection at
the $k$'th step.  Hence, there is positive probability
for no intersections at all.  Given that there is an
intersection at the $k$'th step, the probability for any further
intersections at time $\ge k$ is the same as the probability
to have an intersection at some time $t\geq 1$.
That proves the EIT.

To get the stronger correct estimate for
$$
\P\left[\sum_{j<k} j\alpha_j= \sum_{j<k} j\beta_j\right] 
\,,
$$
we need a bit of Fourier analysis.  The Fourier transform for the
density distribution of $j\alpha_j$ on $\Z$ is
$\cos(jx)$.
Since the Fourier transform transforms convolution into
multiplication, we are interested in the maximum of
$\int_{-\pi}^{\pi}\exp(ixn)\prod_{j<k}\cos(jx)$ over $n\in\Z$.
That is bounded above by
$$
\int_{-\pi}^{\pi}\prod_{j<k}|\cos(jx)|
\,.
$$
Since $\cos(x+\pi)=-\cos(x)$ and $\cos(x)=\cos(-x)$, this is the same as
$$
4\int_0^{\pi/2}\prod_{j<k}|\cos(jx)|
\,.
$$
Let $f(x)$ be $\min\{|x-j\pi|:j\in\Z\}$.
Observe that $|\cos(x)|\leq \exp\left(-cf(x)^2\right)$, for some $c>0$.
So we have to estimate
$$
\int_{0}^{\pi/2}\exp\left(-c\sum_{j<k}f(jx)^2\right)
\,.
$$
We first estimate the integral in the interval $x\in[0,1/k]$:
$$
\int_{0}^{1/k}\exp\left(-c\sum_{j<k}f(jx)^2\right)
=
\int_{0}^{1/k}\exp\left( -c\sum_{j<k}j^2x^2\right)
\leq
\int_{0}^{1/k}\exp\left( -c'k^3x^2\right)
= O\left(k^{-3/2}\right)
\,,
$$
for some constant $c'>0$.
For the other interval $x\in[1/k,\pi/2]$,
note that $\sum_{j<k} f(jx)^2$ is at least $kc''$ for some constant $c''>0$.
Hence the integral
$ \int_{1/k}^{\pi/2}\exp\left(-c\sum_{j<k}f(jx)^2\right)$
is exponentially small in $k$.
\qed

\noindent
{\bf Remarks:}
\begin{enumerate}
\item
The volume of a ball in $G_H$, of radius $n$ is $cn^4$ . $\Z^4$ has the same
volume growth. Yet the probability that
two oriented random walks in $\Z^4$ intersects after $n$ steps is of
order $n^{-3/2}$ and only $n^{-2}$ on $G_H$.

\item
The probability simple random walk on $G_H$ is at
the origin, by time $n$, is of order $cn^{-2}$ (as in $\Z^4$). This can be
used to show that two simple random walk paths will intersect infinitely
often.
\end{enumerate}

\noindent
{\bf Acknowledgements:} Thanks to Ran Raz and Ilya Rips for useful discussions.

\bigskip
\filbreak
\begingroup
{
\small
\begin{sc}
\parindent=0pt\baselineskip=12pt
\def\email#1{\par\qquad {\tt #1} \smallskip}
\def\emailwww#1#2{\par\qquad {\tt #1}\par\qquad {\tt #2}\medskip}

The Weizmann Institute of Science,
Rehovot 76100, Israel
\emailwww{itai@wisdom.weizmann.ac.il}
{http://www.wisdom.weizmann.ac.il/$\sim$itai/}

\emailwww{schramm@wisdom.weizmann.ac.il}
{http://www.wisdom.weizmann.ac.il/$\sim$schramm/}
\end{sc}
}
\filbreak

\end{document}